\documentclass{article}
 \usepackage{amsmath}

\newtheorem{thm}{Theorem}[section]
\newtheorem{cor}[thm]{Corollary}
\newtheorem{prop}[thm]{Proposition}

\newtheorem{lem}[thm]{Lemma}

\newtheorem{ex}[thm]{Example}

\newcommand{\be}{\begin{equation}}
\newcommand{\ee}{\end{equation}}
\newcommand{\ben}{\begin{enumerate}}
\newcommand{\een}{\end{enumerate}}
\newcommand{\beq}{\begin{eqnarray}}
\newcommand{\eeq}{\end{eqnarray}}
\newcommand{\beqn}{\begin{eqnarray*}}
\newcommand{\eeqn}{\end{eqnarray*}}

\newcommand{\pa}{\partial}

\newcommand{\qed}{\hspace*{\fill}Q.E.D.}  

\begin{document}
\title{On a Class of Two-Dimensional Douglas and Projectively Flat Finsler Metrics}
\author{Guojun Yang}
\date{}
\maketitle
\begin{abstract}
In this paper, we study a  class of two-dimensional Finsler
metrics defined by a Riemannian metric $\alpha$ and a $1$-form
$\beta$.   We characterize those metrics which are  Douglasian or
locally projectively flat by some equations. In particular, it
shows that the known fact that $\beta$ is always closed for those
metrics in higher dimensions is no longer true in two dimensional
case. Further, we determine the local structures of
two-dimensional
 $(\alpha,\beta)$-metrics  which are Douglassian, and some
 families
 of examples are given for projectively flat classes with $\beta$
 being not closed.

 \


\end{abstract}

\section{Introduction}

 Projective Finsler geometry studies equivalent Finsler metrics
  on a same manifold with the same geodesics as points (\cite{Dou}). Douglas
  curvature ({\bf D}) is an important projective invariants in projective Finsler geometry. A Finsler metric is called
  {\it Douglasian} if ${\bf D}=0$, and
   {\it locally projectively flat} if at every point, there are local
coordinate systems in which geodesics are straight.  It is known
that a  locally projectively flat Finsle metric can be
characterized by ${\bf D}=0$ and vanishing Weyl curvature.
  As we know, the locally projectively flat class of Riemannian metrics
  is very limited, nothing but the class
  of constant sectional
  curvature (Beltrami Theorem). However, the class of locally projectively flat Finsler metrics is very
  rich. It is known that locally projectively flat Finsler metrics must be Douglassian, but  Douglas metrics  are
  not necessarily locally projectively flat. Therefore, it is
a natural problem to study and classify Finsler metrics which are
Douglasian or locally projectively flat. For this problem, we can
only investigate some special classes of Finsler metrics.

In this paper, we shall consider a special class of Finsler metrics defined by a Riemannian metric
 $\alpha=\sqrt{a_{ij}(x)y^iy^j}$ and a $1$-form $\beta=b_i(x)y^i$ on a manifold $M$. Such metrics
 are called $(\alpha,\beta)$-metrics. An $(\alpha,\beta)$-metric can be expressed in the following form:
 $$F=\alpha \phi(s),\ \ s=\beta/\alpha,$$
where $\phi(s)>0$ is a $C^{\infty}$ function on $(-b_o,b_o)$. It
is known that $F$ is a  regular Finsler metric (defined on the
whole $TM-\{0\}$ and positive definite) for any $(\alpha, \beta)$
with $ \|\beta\|_{\alpha} < b_o$ if and only if
 \be\label{j1}
 \quad \phi(s)-s\phi'(s)+(\rho^2-s^2)\phi''(s)>0,  \quad (|s| \leq \rho <b_o),
 \ee
where $b_o$ is a constant. If $\phi$ does not satisfy (\ref{j1}),
then  $F= \alpha \phi (\beta/\alpha)$ is singular.

Randers metrics are a special class of $(\alpha,\beta)$-metrics.
It is known that a  Randers metric $F= \alpha+\beta$ is a Douglas
metric if and only if $\beta$ is closed (\cite{BaMa}), and it is
locally projectively flat if and only if $\alpha$ is locally
projectively flat and $\beta$ is closed (\cite{BaMa} \cite{CS}).
Usually we call $F=\alpha\phi(s)$ with $\phi(s)=\epsilon
s+\sqrt{1+ks^2}$, where $k,\epsilon$ are constants, a Finsler
metric of {\it Randers type}, which is essentially a Randers
metric.

$(\alpha,\beta)$-metrics are computable and it has been shown that
$(\alpha,\beta)$-metrics have a lot of  special geometric
properties (\cite{KA}--\cite{SY} \cite{Y2}).  In \cite{LSS}
\cite{Shen1}, the authors study and characterize
$(\alpha,\beta)$-metrics which are respectively Douglasian and
locally projectively flat in dimension $n\geq 3$. However, the
two-dimensional case remains open. In this paper, we will solve
this problem in two-dimensional case, and meanwhile give their
local structures in part.

 \begin{thm}\label{th1}
  Let $F=\alpha \phi(s)$, $s=\beta/\alpha$, be a  regular
  $(\alpha,\beta)$-metric on an open subset $U\subset R^2$, where $\phi(0)=1$. Suppose that
    $\beta$ is not parallel with respect to $\alpha$ and $F$ is not
of Randers type. Let $F$ be Douglasian or locally projectively
flat. Then we have one of the following two cases:
 \ben
  \item[{\rm (i)}] $\phi(s)$  satisfies
   \be\label{yg3}
    \big\{1+(k_1+k_3)s^2+k_2s^4\big\}\phi''(s)=(k_1+k_2s^2)\big\{\phi(s)-s\phi'(s)\big\},
   \ee
  where  $k_1,k_2,k_3$ are
constants satisfying
 \be\label{y5200}
 k_2\ne\frac{(2k_1+3k_3)(3k_1+2k_3)}{25},\ \ \ k_2\ne k_1k_3.
 \ee
Further, $\beta$ must be closed.

 \item[{\rm (ii)}] $F$ can be written as
 \be\label{c6}
   F=\widetilde{\alpha}\pm\frac{\widetilde{\beta}^2}{\widetilde{\alpha}},
   \ \ \ \ \ \  \big(\widetilde{\alpha}:=\sqrt{\alpha^2-k\beta^2},\ \
    \widetilde{\beta}:=c\beta\big),
 \ee
    where $k,c$ are constants with  $c\ne 0$. In this case, $\beta$ is generally not closed.
 \een
 \end{thm}

In dimension $n>2$ in Theorem \ref{th1}, it is proved  in
\cite{LSS} \cite{Shen1} that, the metric in Theorem \ref{th1} must
be given by (\ref{yg3}) with $k_2\ne k_1k_3$, and $\beta$ must be
closed. In Theorem \ref{th3} and Theorem \ref{th02} below, we give
general characterizations for
 two-dimensional $(\alpha,\beta)$-metrics (might be singular)
 which are Douglasian and locally projectively
flat respectively.

\

Next we consider the local structure of the Douglas metrics in
Theorem \ref{th1}. By using some deformations on $\alpha$ and
$\beta$, we can determine the local structure of two-dimensional
regular
 Douglas $(\alpha,\beta)$-metrics,
which is shown in the following two theorems. For the local
structure of the singular Douglas classes in Theorem
\ref{th3}(iii) and (iv) below, we will have a discussion in
Section \ref{sec6}.

\begin{thm}\label{th2}
 Let $F=\alpha\pm\beta^2/\alpha$ be a two-dimensional
 regular
 Douglas $(\alpha,\beta)$-metric with $\beta\ne 0$. Then
 $\alpha$ and $\beta$ can be locally written as
  \beq
 \alpha^2&=&\frac{1}{(1\mp B)^3}\Big\{\frac{B}{u^2+v^2}\big[(y^1)^2+(y^2)^2\big]\mp 9(1\pm B+B^2)\beta^2\Big\},\label{yc1}\\
 \beta&=&\frac{B}{(1\pm 2B)^{\frac{3}{2}}}\frac{uy^1+vy^2}{u^2+v^2},\label{yc2}
  \eeq
  where
   $u=u(x),v=v(x)$ are  scalar functions such that
    $$f(z)=u+iv, \ \ z=x^1+ix^2$$
    is a complex analytic function, and $B=B(x)$ is a scalar
    function satisfying $0<B<1$ if $F=\alpha+\beta^2/\alpha$ and $0<B<1/2$ if
    $F=\alpha-\beta^2/\alpha$.
\end{thm}

 We see that in Theorem \ref{th2}, the metric is
determined by the triple parametric scalar functions $(B,u,v)$,
where $u$ and $v$ are a pair of complex conjugate functions. We
will prove Theorem \ref{th2} by using  Corollary \ref{cor42} and
the result in \cite{Y} (also see \cite{Y1}).

\begin{thm}\label{th001}
Let $F=\alpha \phi(s)$, $s=\beta/\alpha$ be a two-dimensional
regular Douglas
  $(\alpha,\beta)$-metric with $\beta\ne 0$, where $\phi(s)$ satisfies (\ref{yg3}) and
 (\ref{y5200}),  and $\phi(0)=1$. Then $\alpha$ and $\beta$ can be
  locally written as
   \be\label{yg79}
 \alpha=\frac{\sqrt{B}}{c}\sqrt{\frac{(y^1)^2+(y^2)^2}{u^2+v^2}},\ \
 \beta=\frac{B(uy^1+vy^2)}{c(u^2+v^2)}, \ \ \
 \big(c:=e^{\int^{B}_0\frac{1}{2}\frac{k_3+k_2t}{1+(k_1+k_3)t+k_2t^2}dt}\big),
   \ee
   where
    $B=B(x)>0,u=u(x),v=v(x)$ are some scalar functions
   which satisfy the following PDEs:
   \be\label{yg81}
  u_1=v_2,\ \ u_2=-v_1,\ \ v_1+\sigma_1v=u\sigma_2,
   \ee
where $u_i:=u_{x^i},v_i:=v_{x^i}$ and $\sigma_i:=\sigma_{x^i}$,
and $\sigma$ is defined by
 \be\label{sigma}
 e^{2\sigma}:=\frac{B}{c^2(u^2+v^2)}.
 \ee
\end{thm}

 We will
prove Theorem \ref{th001} by (\ref{yg4}) and  the result in
\cite{Y} (also see \cite{Y1}).  The metric in Theorem \ref{th001}
is determined by the triple parametric scalar functions $(B,u,v)$
which satisfy (\ref{yg81}). It seems hard to obtain the complete
solutions of the PDEs (\ref{yg81}). However, we can give some
special solutions of (\ref{yg81}). For example, the following
triple is a solution
 $$
\sigma=x^1,\ \ u=(c_2\sin x^2-c_1\cos x^2)e^{-x^1},\ \  v=(c_1\sin
x^2+c_2\cos x^2)e^{-x^1},
 $$
where $c_1,c_2$ are constants, and then $B$ is determined by
(\ref{sigma}).

\

Now we consider the local structure of the locally projectively
flat metrics in Theorem \ref{th1}. The local structure of $F$
determined by (\ref{yg3}) with $k_2\ne k_1k_3$ (for the dimension
$n\ge 2$) has been solved in \cite{Yu} (also see another way in
\cite{Y2}). However, it seems difficult to determine the local
structure of $F$ determined by (\ref{c6}). By using (\ref{y0006})
and (\ref{y00060}) with $\tau=0$, we can construct the following
example, which can be directly verified. We omit the details.

\begin{ex}
 Let $F=\alpha\pm\beta^2/\alpha$ be a two-dimensional
 $(\alpha,\beta)$-metric. Suppose $\alpha$ and $\beta$ take the
 form
\be\label{yg95}
  \alpha=e^{\sigma(x)}\sqrt{(y^1)^2+(y^2)^2}, \ \
  \beta=e^{\sigma(x)}(\xi(x) y^1+\eta(x) y^2).
  \ee
 where $\xi=\xi(x),\eta=\eta(x),\sigma=\sigma(x)$ are some scalar functions. Define
  \beq
 \xi(x)&=&\pm\frac{1}{\sqrt{2}}\frac{x^2+c_2}{\sqrt{c_3\mp(x^1+c_1)^2\mp(x^2+c_2)^2}},\label{y131}\\
 \eta(x)&=&\mp\frac{1}{\sqrt{2}}\frac{x^1+c_1}{\sqrt{c_3\mp(x^1+c_1)^2\mp(x^2+c_2)^2}},\label{y132}\\
 \sigma(x)&=&\ln\big[c_3\mp(x^1+c_1)^2\mp(x^2+c_2)^2\big]+c_4,\label{y133}
  \eeq
where $c_1,c_2,c_3$ are constants with $c_3>0$. Then $F$ is
projectively flat with $\beta$ being not closed.
\end{ex}

For the singular projectively flat classes in Theorem
\ref{th02}(iii) and (iv) below, we also construct some examples
with $\beta$ being not closed (see Example \ref{ex01} and
\ref{ex02} below). As we have shown, it seems an obstacle to
determine the local structure of the projectively flat classes
when $\beta$ is not closed.

\

\noindent{\bf Open Problem:} Determine the local structure of a
two-dimensional
 $(\alpha,\beta)$-metric $F=\alpha\pm\beta^2/\alpha$ which is locally projectively
 flat.

\section{Preliminaries}

Let $F=F(x,y)$ be a Finsler metric on an $n$-dimensional manifold $M$.
      The geodesic coefficients are defined by
 \beq \label{G1}
 G^i:=\frac{1}{4}g^{il}\big \{[F^2]_{x^ky^l}y^k-[F^2]_{x^l}\big \}.
 \eeq

    A Finsler  metric $F=F(x,y)$ is called a {\it
Douglas metric} if the spray coefficients $G^i$ are in the
following form:
 \beq\label{G2}
   G^i=\frac{1}{2}\Gamma_{jk}^i(x)y^jy^k+P(x,y)y^i,
 \eeq
where $\Gamma_{jk}^i(x)$ are local functions on $M$ and $P(x,y)$
is a local positively homogeneous function of degree one in $y$.
It is easy to see that $F$
 is a Douglas metric if and only if $G^iy^j-G^jy^i$ is a
 homogeneous polynomial in $(y^i)$ of degree three, which by
 (\ref{G2}) can be written as (\cite{BaMa}),
  $$G^iy^j-G^jy^i=\frac{1}{2}(\Gamma^i_{kl}y^j-\Gamma^j_{kl}y^i)y^ky^l.$$

According to G. Hamel's result, a Finsler metric $F$ is
projectively flat in $U$ if and only if
 $$F_{x^my^l}y^m-F_{x^l}=0.$$
 The above formula implies that $G^i=Py^i$ with $P$ given by
  $$P=\frac{F_{x^m}y^m}{2F}.$$

For a Riemannian metric $\alpha =\sqrt{a_{ij}y^iy^j}$ and a
$1$-form $\beta=b_iy^i$ on a manifold $M$, let $\nabla \beta =
b_{i|j} y^i dx^j$  denote the covariant derivatives of $\beta$
with respect to $\alpha$. Put
 $$r_{ij}:=\frac{1}{2}(b_{i|j}+b_{j|i}),\ \ s_{ij}:=\frac{1}{2}(b_{i|j}-b_{j|i}),\ \
 r_j:=b^ir_{ij},\ \ s_j:=b^is_{ij},\ \ s^i:=a^{ik}s_k,$$
 where $b^i:=a^{ij}b_j$ and $(a^{ij})$ is the inverse of
 $(a_{ij})$.

Consider an $(\alpha,\beta)$-metric $F =\alpha \phi
(\beta/\alpha)$.
 By (\ref{G1}), the spray coefficients $G^i$ of $F$
are given by (\cite{CS} \cite{KA} \cite{Ma} \cite{Shen1} \cite{Shen2}):
  \be\label{y20}
  G^i=G^i_{\alpha}+\alpha Q s^i_0+\alpha^{-1}\Theta (-2\alpha Q
  s_0+r_{00})y^i+\Psi (-2\alpha Q s_0+r_{00})b^i,
  \ee
where $s^i_j=a^{ik}s_{kj}, s^i_0=s^i_ky^k,
s_i=b^ks_{ki},s_0=s_iy^i$, and
 $$
  Q:=\frac{\phi'}{\phi-s\phi'},\ \
  \Theta:=\frac{Q-sQ'}{2\Delta},\ \
  \Psi:=\frac{Q'}{2\Delta},\ \ \Delta:=1+sQ+(b^2-s^2)Q'.
 $$

By (\ref{y20}) one  can see that $F=\alpha\phi(\beta/\alpha)$ is a
Douglas metric if and only if
 \be\label{y21}
 \alpha Q (s^i_0y^j-s^j_0y^i)+\Psi (-2\alpha
 Qs_0+r_{00})(b^iy^j-b^jy^i)=\frac{1}{2}(G^i_{kl}y^j-G^j_{kl}y^i)y^ky^l,
\ee
  where $G^i_{kl}:=\Gamma^i_{kl}-\gamma^i_{kl}$,  $\Gamma^i_{kl} $ are given in (\ref{G2}) and $
  \gamma^i_{kl}:=\pa^2G^i_{\alpha}/\pa y^k\pa y^l.$

Further,  $F=\alpha\phi(\beta/\alpha)$ is projectively flat on $U\subset R^n$ if and only if
\be
(a_{ml}\alpha^2-y_my_l)G^m_{\alpha}+\alpha^3Qs_{l0}+\Psi\alpha(-2\alpha
 Qs_0+r_{00})(\alpha b_l-sy_l)=0,\label{y21*}
\ee
where $y_l=a_{ml}y^m$.

\section{Equations in a Special Coordinate System}\label{s3}
In order to prove Theorems \ref{th3} and \ref{th02} below, one has
to simplify  (\ref{y21}) and (\ref{y21*}). The main technique is
to fix a point and choose a   special  coordinate
 system $(s,y^a)$  as  in \cite{Shen1} \cite{Shen2}.

Fix an arbitrary point $x\in M$ and take  an orthonormal basis
  $\{e_i\}$ at $x$ such that
   $$\alpha=\sqrt{\sum_{i=1}^n(y^i)^2},\ \ \beta=by^1.$$
Then we change coordinates $(y^i)$ to $(s, y^a)$ such that
  $$\alpha=\frac{b}{\sqrt{b^2-s^2}}\bar{\alpha},\ \
  \beta=\frac{bs}{\sqrt{b^2-s^2}}\bar{\alpha}, $$
where $\bar{\alpha}=\sqrt{\sum_{a=2}^n(y^a)^2}$. Let
 $$\bar{r}_{10}:=r_{1a}y^a, \ \ \bar{r}_{00}:=r_{ab}y^ay^b, \ \
 \bar{s}_0:=s_ay^a.$$
We have $\bar{s}_0=b\bar{s}_{10},s_1=bs_{11}=0$. The following lemmas are trivial.

\begin{lem}\label{lem5.1}
In the special local coordinate system at $x$ as mentioned above, if
 $b=constant$, then $r_{11}=0, r_{1a}+s_{1a}=0$ at $x$.
\end{lem}

\begin{lem}{\rm(\cite{Shen2})}
  For $n\ge 2$, suppose $p+q\bar{\alpha}=0$, where
  $p=p(\bar{y})$ and $q=q(\bar{y})$ are homogeneous
  polynomials in $\bar{y}=(y^a)$, then $p=0,q=0$.
 \end{lem}

By
  \cite{LSS} and \cite{Shen1} we have the following two propositions.

\begin{prop}\label{prop1}
  $(n=2)$ An $(\alpha,\beta)$-metric
 $F=\alpha\phi(\beta/\alpha)$ is a Douglas metric if and only if at each point $x$,
 there is a suitable coordinate system such that at $x$,
 there exist numbers  $G^i_{jk}$ $(i,j,k=1,2)$ which are independent of $s$ such that
 \be\label{y27}
  \frac{s^2}{2(b^2-s^2)}(G^1_{11}-G^2_{12}-G^2_{21})+\frac{1}{2}G^1_{22}=
   b\Psi (\frac{s^2}{b^2-s^2}r_{11}+r_{22}),
 \ee
 \be\label{y28}
 \frac{1}{b^2-s^2}\big[2\Psi
 (b^2-s^2)-1\big]b^3Qs_{12}-2b\Psi r_{12}s=\frac{G^2_{11}}{2(b^2-s^2)}s^3+\frac{1}{2}(G^2_{22}-G^1_{12}-G^1_{21})s.
 \ee
 \end{prop}

 \begin{prop}\label{prop2}
  $(n=2)$ An $(\alpha,\beta)$-metric
 $F=\alpha\phi(\beta/\alpha)$ is projectively flat if and only if
 \be\label{g26}
  \frac{s^2}{2(b^2-s^2)}(-\widetilde{G}^1_{11}+2\widetilde{G}^2_{12})-\frac{1}{2}\widetilde{G}^1_{22}=
   b\Psi (\frac{s^2}{b^2-s^2}r_{11}+r_{22}),
 \ee
 \be\label{g27}
 \frac{1}{b^2-s^2}\big[2\Psi
 (b^2-s^2)-1\big]b^3Qs_{12}-2b\Psi r_{12}s=-\frac{\widetilde{G}^2_{11}}{2(b^2-s^2)}s^3+
 \frac{1}{2}(-\widetilde{G}^2_{22}+2\widetilde{G}^1_{12})s,
 \ee
 where $\widetilde{G}^i_{jk}:=\frac{\pa^2G_{\alpha}^i}{\pa y^j\pa y^k}$ are
 the connection coefficients of $\alpha$.
 \end{prop}

 Comparing (\ref{y27}) and (\ref{g26}), (\ref{y28}) and (\ref{g27}), it is easy to see that if
$G^i_{jk}=G^i_{kj}$, then $\widetilde{G}^i_{jk}=-G^i_{jk}$. So if
we can solve $G^i_{jk}$ from (\ref{y27}) and (\ref{y28}), then we
can solve $\widetilde{G}^i_{jk}$ from (\ref{g26}) and (\ref{g27}).
In the following we only consider (\ref{y27}) and (\ref{y28}),
from which we will solve $G^i_{jk}$.

\section{Douglas $(\alpha,\beta)$-metrics}

In this section, we characterize two-dimensional
$(\alpha,\beta)$-metrics (might be singular) which are Douglas
metrics. We have the following theorem.

\begin{thm}\label{th3}
  Let $F=\alpha \phi(s)$, $s=\beta/\alpha$, be an
  $(\alpha,\beta)$-metric on an open subset $U\subset R^2$, where $\phi(0)=1$. Suppose that
    $\beta$ is not parallel with respect to $\alpha$ and $F$ is not
of Randers type. Then $F$ is  a Douglas metric on $U$ if and only
if $F$ lies in one of the following four classes:
 \ben
  \item[{\rm (i)}] $\phi(s)$  satisfy (\ref{yg3}) with $k_2\ne k_1k_3$ and $\beta$
  satisfies
   \be\label{yg4}
     b_{i|j}=2\tau
     \big\{(1+k_1b^2)a_{ij}+(k_2b^2+k_3)b_ib_j\big\},
   \ee

 \item[{\rm (ii)}] $\phi(s)$ and $\beta$ satisfy
\beq
     \phi(s)&=&\sqrt{1-ks^2} + \frac{cs^2}{\sqrt{1-k s^2}},\label{y5}\\
     r_{ij}&=&2\tau
     \Big\{[1+ (2c-k)b^2]a_{ij}-\big[k+3c -(k+c)k b^2\big]b_ib_j\Big\}\nonumber\\
     &&+d(b_is_j+b_js_i),\label{y6}
    \eeq
   where $\tau=\tau(x)$ is a scalar function,
   $k, c$ are constants with $c \not=0$ and $1-k b^2 \geq 0$,  and $d=d(x)$ is given by
   \be\label{y7}
    d=\frac{ 3c-k - (2c-k)k b^2}{1-(k+c)b^2}.
   \ee

 \item[{\rm (iii)}] ($b= constant$)  $\phi(s)$ and $\beta$ satisfy
 \beq
 \phi(s)&=& \frac{\sqrt{b^2-s^2}}{b} + \sqrt{b^2-s^2} \int_0^s
\frac{c}{(b^2-t^2)^{3/2} } \Big ( \frac{t^2}{1-k t^2} \Big )^m
dt,\label{y8}\\
  r_{ij}&=&-\frac{1}{b^2}(b_is_j+b_js_i),\label{y9}
  \eeq
  where $c, k$ are constants and $m\ge 1$ is an integer.

 \item[{\rm (iv)}] ($b=constant$)  $\phi(s)$ and $\beta$ satisfy
 \beq
\phi(s)&=& \frac{\sqrt{b^2-s^2}}{b} + \sqrt{b^2-s^2}\int_0^s
\frac{c}{(b^2-t^2)^{3/2} } \Big ( \frac{t^2}{ 1-kt^2 } \Big
)^{m-1/2} dt,\label{y10}\\
 r_{ij}&=&-\frac{1}{b^2}(b_is_j+b_js_i),\label{y11}
 \eeq
 where $c, k$ are constants and $m\ge 1$ is an integer.
 \een
 \end{thm}

By Theorem \ref{th3}(ii), we can easily get the following
corollary.

\begin{cor}\label{cor42}
 Let $F=\alpha\pm \beta^2/\alpha$ be  two-dimensional
 $(\alpha,\beta)$-metric. Then $F$ is a Douglas metric if and only
 if  $\beta$ satisfies
 \be\label{y0006}
  r_{ij}=2\tau\big\{(1\pm 2b^2)a_{ij}\mp
  3b_ib_j\big\}
  +\frac{3}{\pm1-b^2}(b_is_j+b_js_i),
 \ee
 where $\tau=\tau(x)$ is a scalar function. Note that $F=\alpha+
 \beta^2/\alpha$ is regular if and only if $b<1$; $F=\alpha-
 \beta^2/\alpha$ is regular if and only if $b<1/2$.
\end{cor}

We prove Theorem \ref{th3} using Proposition \ref{prop1}. The
proof can be divided  into two cases $(r_{11}, r_{22}) \not=(0,
0)$ and $(r_{11}, r_{22}) = (0, 0)$.

\subsection{ $(r_{11},r_{22})\ne (0,0)$}\label{sec4}

In this case, we wil obtain two classes:  Theorem \ref{th3} (i)
and Theorem \ref{th3} (ii).

First,
  (\ref{y27}) can be written in the following form
 \be\label{y29}
  2\Psi=\frac{\lambda s^2+\mu (b^2-s^2)}{\delta s^2+\eta
  (b^2-s^2)},
 \ee
where $\lambda,\mu,\delta,\eta$ are numbers independent of $s$. By
(\ref{y27}) and (\ref{y29}), it is easy to prove  that if
$\lambda\eta-\mu\delta\ne 0$, then for some scalar $\tau=\tau(x)$,
we have (see also \cite{LSS})
 \be\label{y30}
 r_{11}=2b^2\delta \tau,\ \ r_{22}=2b^2\eta \tau.
 \ee
One can see that if an $(\alpha,\beta)$-metric
$F=\alpha\phi(\beta/\alpha)$ is not of Randers type, then
$\lambda\eta-\mu\delta\ne 0$.

Now we put
 $$ a_0:=1,\ \ a_i:=\frac{\phi^{(i)}(0)}{i!}, \ \ i=1,2,\cdots.$$
By (\ref{y29}), it has been proved in \cite{Shen1} that if
$2a_4+a_2^2= 0$, then $F$ is of Randers type.
Thus we may assume that $2a_4+a_2^2\ne
0$. Then there is a  scalar $\epsilon=\epsilon(x)\ne 0$ such that
 \be\label{y32}
 \mu=k_1\epsilon,\ \ \eta=(1+k_1b^2)\epsilon,\ \
 \lambda=(k_1+k_2b^2)\epsilon,\ \
 \delta=(1+(k_1+k_3)b^2+k_2b^4)\epsilon,
\ee
 \be\label{y33}
  2\Psi=\frac{k_1+k_2s^2}{1+k_1b^2+(k_3+k_2b^2)s^2},
 \ee
 where $k_1,k_2,k_3$ are some constants determined by
 \be\label{y34}
  k_1=2a_2,\ \
  k_2=\frac{2(a_4a_2^2-5a_2a_6+12a_4^2)}{2a_4+a_2^2},
   \ \ k_3=-\frac{11a_2a_4+5a_6+3a_2^3}{2a_4+a_2^2}.
 \ee
Note that $k_2-k_1k_3\ne 0$ is equivalent to $2a_4+a_2^2\ne 0$.
Since $F$ is not of Randers type, we get $k_2-k_1k_3\ne 0$.

Plugging (\ref{y33}) into  (\ref{y28}),  we get
 \beq
  &&-2bs(b^2-s^2)(k_1+k_2s^2)r_{12}-2b^3(1+k_1s^2+k_3s^2+k_2s^4)Qs_{12}+
  \nonumber\\
&& \qquad \qquad
s(1+k_1b^2+k_3s^2+k_2b^2s^2)\big\{(b^2-s^2)\xi-G^2_{11}s^2\big\}=0,\label{y36}
 \eeq
 where $\xi:=G^1_{12}+G^1_{21}-G^2_{22}$.
Now plug the Taylor expansion of $\phi(s)$ into (\ref{y36}) and
let $p_i$ be the coefficients of $s^i$ in (\ref{y36}). By
$p_1=0,p_3=0,p_5=0$ we have the following cases.

   (i) If
   \be\label{yg}
   1+(k_1+k_3)b^2+k_2b^4\ne 0,
   \ee
  then
  \be\label{y41}
   r_{12}=\frac{b^2}{15}\frac{k_1(3k_1^2+2k_1k_3+10k_2)b^4+(21k_1^2+14k_1k_3-5k_2)b^2+18k_1-3k_3}
   {1+(k_1+k_3)b^2+k_2b^4}s_{12},
  \ee
  \be\label{y42}
   G^2_{11}=\frac{2b}{15}\frac{(3k_1^2k_3+2k_1k_3^2-2k_2k_3-18k_1k_2)b^4-5(3k_1^2+k_2+2k_1k_3)b^2-15k_1}
   {1+(k_1+k_3)b^2+k_2b^4}s_{12},
  \ee
  \be\label{y43}
   \xi=\frac{2bk_1}{15}\frac{(3k_1^2+2k_1k_3+10k_2)b^4+6(3k_1+2k_3)b^2+15}{1+(k_1+k_3)b^2+k_2b^4}s_{12}.
  \ee

  (ii) If
  \be\label{y44}
   1+(k_1+k_3)b^2+k_2b^4=0,
  \ee
  then $b=constant$.
By Lemma \ref{lem5.1} we get $r_{12}+s_{12}=0$.
  If $s_{12}\ne 0$, then we get
   \be\label{y48}
    k_2=\frac{3+k_1b^2}{2b^4},\ \ k_3=-\frac{5+3k_1b^2}{2b^2}.
   \ee

\bigskip

\noindent{\bf Case 1}:  $s_{12}=0$. This implies Theorem \ref{th3}
(i).

\bigskip

\noindent
If (\ref{yg}) holds, then $r_{12}=0$ by
 (\ref{y41}). If (\ref{y44}) holds, we also get $r_{12}=0$. In
 both cases, we have $G^2_{11}=\xi=0$ by plugging
 $r_{12}=s_{12}=0$ into (\ref{y28}). Thus
 (\ref{y28}) becomes trivial.
 By $r_{12}=0$,
 (\ref{y30}) and (\ref{y32}), we get
 the expression of $b_{i|j}$ in (\ref{yg4}). Further,
 (\ref{y33}) can be written in the form (\ref{yg3}) with $k_2\ne k_1k_3$.
 This class belongs to Theorem \ref{th3}(i).

\bigskip
\noindent{\bf Case 2}:
 $s_{12}\ne 0$. This implies Theorem \ref{th3} (ii).

\noindent {\bf Case 2A}.
Assume that  (\ref{yg}) holds.  We plug (\ref{y41}), (\ref{y42}) and (\ref{y43}) into
  (\ref{y36}), and then we obtain
  \be\label{y49}
   Q=-\frac{1}{15}\frac{\big\{(3k_1^2k_3+2k_1k_3^2-2k_2k_3-18k_1k_2)s^4-5(3k_1^2+k_2+2k_1k_3)s^2-15k_1)\big\}s}
    {1+(k_1+k_3)s^2+k_2s^4}.
  \ee
By (\ref{y49}) and  (\ref{y33}) we have
 \be\label{y52}
 k_2=\frac{(2k_1+3k_3)(3k_1+2k_3)}{25}.
 \ee
 Plug (\ref{y52}) into (\ref{y49}) and we get
 \be\label{y54}
  \phi(s)=\frac{1}{\sqrt{5}}\frac{5+(4k_1+k_3)s^2}{\sqrt{5+(3k_1+2k_3)s^2}},
  \ee
where $k_1\ne k_3$ since (\ref{y52}) and $k_2\ne k_1k_3$. Letting
$ k_1 = 2c-k$ and $k_3 = -3c-k$ in (\ref{y54}),  we get
(\ref{y5}).
 Substituting
(\ref{y52}) into (\ref{y41}) gives
 \be\label{y55}
 r_{12}=\frac{b^2\big[k_1(3k_1+2k_3)b^2+6k_1-k_3\big]}{5+(2k_1+3k_3)b^2}s_{12}.
 \ee
Letting $k_1 = 2c-k$ and $k_3 = -3c -k$ and using (\ref{y30}),
 we obtain (\ref{y6}).

\bigskip
\noindent {\bf Case 2B}. Assume that  (\ref{y44}) holds. Then
$r_{12}=-s_{12}$ and (\ref{y48}) holds. It is easy for us to get
 \beq
     \phi(s)&=&\frac{b+\tilde{c}s^2}{\sqrt{b^2-s^2}},\label{y12}\\
     r_{ij}&=&2\tilde{\tau}(b^2a_{ij}-b_ib_j)-\frac{1}{b^2}(b_is_j+b_js_i).\label{y13}
 \eeq
This class is a special case of Theorem \ref{th1}(ii).

\subsection{ $(r_{11},r_{22})= (0,0)$}\label{sec5}

 Since $\beta$ is not parallel and $(r_{11},r_{22})= (0,0)$, we will see that $s_{12}\ne 0$ from
 the following proof to different cases.
It follows from   (\ref{y27}) that
\[ G^1_{22}=0, \ \ \ \ \ \ G^1_{11}=G^2_{12}+G^2_{21}.\]
Plugging the expressions of $Q$ and $\Psi$ into (\ref{y28}) yields
 \beq
 &&s(b^2-s^2)\big[2br_{12}+G^2_{11}s^2-(b^2-s^2)\xi\big]\phi''\nonumber\\
 &&\ \ \ \ +
 s\big[G^2_{11}s^2-(b^2-s^2)\xi\big](\phi-s\phi')+2b^3s_{12}\phi'=0,\label{yy49}
 \eeq
where $\xi:=G^1_{12}+G^1_{21}-G^2_{22}$.

Let
$$\phi=a_0+\sum^{h}_{i=1}a_is^i+o(s^h),\ \ a_0=1,$$
where $h$ is a sufficiently large integer.
Plugging  the above Taylor series into (\ref{yy49}) we obtain a power series $\sum_k p_k s^k =0$.
It is
easily seen that the coefficient $p_k$ of $s^k$  is given by
 \be\label{yy50}
  p_k=A_1r_{12}+A_2G^2_{11}+A_3\xi+A_4s_{12},
 \ee
 where
  \beqn
  A_1:&=&b\big[(k-1)(k-2)a_{k-1}-k(k+1)a_{k+1}b^2\big],\\
  A_2:&=&\frac{1}{2}(k-2)\big[(k-4)a_{k-3}-(k-1)a_{k-1}b^2\big],\\
  A_3:&=&\frac{1}{2}\big[k(k+1)a_{k+1}b^4-(2k-1)(k-2)a_{k-1}b^2+(k-2)(k-4)a_{k-3}\big],\\
  A_4:&=&-(k+1)a_{k+1}b^3,
  \eeqn
 and $a_i=0$ if $i<0$. In particular, we have
  $$
 p_0=-a_1b^3s_{12}.
  $$
  So if $s_{12}\ne 0$, then by $p_0=0$ we have $a_1=0$.

\bigskip

\noindent
 {\bf Case I:} Suppose $db\ne 0$. We will prove that one case
 belongs to Theorem \ref{th3} (ii) with the
 scalar $\tau=\tau(x) =0$, and other cases are excluded.

 Solving the system $p_1=0,p_3=0,p_5=0$ yields the following
 three cases:

 (i) If $2a_4\ne -a_2^2$, then we get (\ref{y41}), (\ref{y42}) and
 (\ref{y43}) by using (\ref{y34}).

  (ii) If $2a_4= -a_2^2$ and
  $2a_6\ne a_2^3$, then
  \be\label{yy56}
   r_{12}=\frac{1}{5}(8a_2b^2-1)s_{12},\ \ G^2_{11}=-4a_2bs_{12},\
   \ \xi=\frac{16}{5}a_2bs_{12}.
  \ee

   (iii) If $2a_4= -a_2^2$ and
  $2a_6=a_2^3$, then
  \be\label{yy57}
   G^2_{11}=-4a_2bs_{12},\ \
   \xi=\frac{4a_2b}{1+2a_2b^2}(r_{12}+s_{12}).
  \ee
It follows from (\ref{y41}) or (\ref{yy56}) that $r_{12}=0$ if
$s_{12}=0$. If (\ref{yy57}) holds and $s_{12}=0$, then we have
$G^2_{11}=0$, and thus we have $r_{12}=0$ by (\ref{y28}) since $F$
is not of Randers type (also see the proof in \cite{LSS}).
Therefore in this case we have $s_{12}\ne 0$.

\bigskip

 \noindent
 {\bf Case IA.}  Suppose $2a_4= -a_2^2$ and
  $2a_6\ne a_2^3$. Then plug (\ref{yy56}) into
  (\ref{yy49}) and by using $db\ne 0$ we get
 $\phi(s)=\sqrt{1+2a_2s^2}$. This case is excluded.

\bigskip
\noindent
 {\bf Case IB.}  Suppose $a_4\ne -\frac{1}{2}a_2^2$.
 Plugging (\ref{y41}), (\ref{y42}) and (\ref{y43}) into (\ref{yy49})
 and by using $db\ne 0$ and $s_{12}\ne 0$ we obtain three ODEs on
 $\phi(s)$, whose discussion of solutions can be divided into the
 following cases.

If $k_3\ne k_1$, then in a similar way as  in section \ref{sec4},
we can easily show that this
 class belongs to  Theorem \ref{th1}(ii) with $\tau=\tau(x)=0$.

 If $k_3=k_1\ne 0$, then  we obtain
 $k_2=k_1^2$, which is
 impossible since $k_2\ne k_1k_3$.

 If $k_1=k_3=0$, then
 it is easy to see that $F$ is of Randers type, which is excluded.

\bigskip
\noindent {\bf Case IC.} Suppose $a_4= -\frac{1}{2}a_2^2$ and
  $a_6=\frac{1}{2}a_2^3$. Then we have (\ref{yy57}). Note that we have $a_1=0$
  since $s_{12}\ne 0$. We will show
  that this case is excluded.

  For the function
  $f(s)=\sqrt{1+2a_2s^2}$, its Taylor
  coefficients  $c_i$ of  $s^i$ ($i\ge0$) are given by
  $$c_{2i+1}=0,\ \ c_{2i}=C^i_{\frac{1}{2}}(2a_2)^i,$$
 where $C^i_{\mu}$ are the
 generalized combination coefficients. So in all $a_{2i+1}$'s
 or $a_{2i}$'s there exist some minimal $m$ such that
 \be\label{y00}
  a_{2m+1}\ne 0,\ \ (m\ge1); \ \ {\text or} \ \
 a_{2m}\ne C^m_{\frac{1}{2}}(2a_2)^m,\ \ (m\ge4).
 \ee

\noindent {\bf Case IC(1).} Assume $a_{2m+1}\ne 0$ in (\ref{y00}).
Then plugging (\ref{yy57}), $a_{2m-3}=0$ and $a_{2m-1}=0$ into
$p_{2m}=0$ (see (\ref{yy50})) yields
 \be\label{yy68}
  a_{1+2m}\big[-2mr_{12}+(4ma_2b^2-2a_2b^2-1)s_{12}\big]=0.
 \ee
 Therefore, it follows from (\ref{yy68})
  that we have
  \be\label{yy70}
  r_{12}=\frac{4ma_2b^2-2a_2b^2-1}{2m}s_{12},
  \ee

\noindent {\bf Case IC(2).} Assume $a_{2m}\ne
C^m_{\frac{1}{2}}(2a_2)^m$ in (\ref{y00}). Plugging (\ref{yy57})
and
$$a_{2m-4}=C^{m-2}_{\frac{1}{2}}(2a_2)^{m-2},\ \
a_{2m-2}=C^{m-1}_{\frac{1}{2}}(2a_2)^{m-1}$$ into $p_{2m-1}=0$
(see (\ref{yy50})) yields
 \be\label{yy69}
  \big[a_{2m}-C^m_{\frac{1}{2}}(2a_2)^m\big]\big[(1-2m)r_{12}+(4ma_2b^2-4a_2b^2-1)s_{12}\big]=0.
 \ee
 Therefore it follows from (\ref{yy69})
  that we have
  \be\label{yy71}
  r_{12}=\frac{4ma_2b^2-4a_2b^2-1}{2m-1}s_{12}.
  \ee

 Finally, plugging (\ref{yy57}) and (\ref{yy70}) or (\ref{yy71})  into
(\ref{yy49}) and using  $db\ne 0$ we get
  $\phi(s)=\sqrt{1+2a_2s^2}$. Thus both cases are excluded.

\bigskip
\noindent
 {\bf Case II:} Suppose $db=0$. We will obtain  Theorem \ref{th3} (iii) and Theorem \ref{th1} (iv).

  By Lemma \ref{lem5.1} we have
 \be\label{yy76}
 r_{12}=-s_{12}.
 \ee
 For the function $f(s)=\frac{1}{b}\sqrt{b^2-s^2}$, its Taylor
  coefficients  $c_i$ of  $s^i$ ($i\ge0$) are given by
  $$c_{2i+1}=0,\ \ c_{2i}=C^i_{\frac{1}{2}}(-\frac{1}{b^2})^i,$$
Since $a_1=0$ and $F$ is not of Randers
 type, in all $a_{2i+1}$'s or $a_{2i}$'s there exist some minimal $m$ such that
 \be\label{y01}
  a_{2m+1}\ne 0,\ \ (m\ge1); \ \ {\text or} \ \
 a_{2m}\ne C^m_{\frac{1}{2}}(-\frac{1}{b^2})^m,\ \ (m\ge1).
 \ee

 \noindent{\bf Case II(1):} Assume $a_{2m+1}\ne 0$ in (\ref{y01}). Plugging
(\ref{yy76}), $a_{2m-3}=0$ and $a_{2m-1}=0$ into $p_{2m}=0$ (see
(\ref{yy50})) yields
 \be\label{yy80}
 \xi=-\frac{2m-1}{mb}s_{12}.
 \ee
 Plugging (\ref{yy80}) and $a_{2m-1}=0$ into $p_{2m+2}=0$  yields
 \be\label{yy81}
 G^2_{11}=\frac{1}{m(2m+1)}\Big[\frac{(2m+3)ba_{2m+3}}{ma_{2m+1}}+\frac{4m^2-3}{b}\Big]s_{12}.
 \ee
Now plug (\ref{yy76}), (\ref{yy80}), (\ref{yy81}) into
(\ref{yy49}), and then we obtain
 \be
\frac{ \phi-s\phi'+(b^2-s^2)\phi''}{ s\phi +(b^2-s^2)\phi'} =
\frac{2m}{s(1-ks^2)},\label{y8****}
 \ee
  where $k$ is a constant determined by $a_{2m+1}$ and $a_{2m+3}$.
Let
\[ \Phi:= s \phi(s)+(b^2-s^2)\phi'(s).\]
Then (\ref{y8****}) becomes
\[ \frac{\Phi'}{\Phi} = \frac{2m}{s(1-ks^2)}.\]
We get
\[ \Phi = c \Big ( \frac{s^2}{1-ks^2} \Big )^m,\]
where $c$ is a constant. Then we can easily get
 \be \phi =
\sqrt{b^2-s^2} \int \frac{c}{(b^2-s^2)^{3/2}} \Big (
\frac{s^2}{1-ks^2} \Big )^m ds.
 \ee
By assumption, $\phi(0)=1$, we get (\ref{y8}). Further, since
$r_{11}=0,r_{22}=0$ and $r_{12}=-s_{12}$, we get (\ref{y9}). This
class belongs to  Theorem \ref{th3}(iii).

\bigskip

\noindent{\bf Case II(2):}  Assume $a_{2m}\ne
C^m_{\frac{1}{2}}(-\frac{1}{b^2})^m$ in (\ref{y01}). Plugging
(\ref{yy76}) and the expressions of $a_{2m-4}$ and $a_{2m-2}$
 into
$p_{2m-1}=0$ yields
 \be\label{yy83}
  \xi=-\frac{4(m-1)}{(2m-1)b}s_{12}.
 \ee
Plugging (\ref{yy83}) and the expressions of $a_{2m-2}$ into
$p_{2m+1}=0$ (see (\ref{yy50})) yields
 \be\label{yy84}
 G^2_{11}=\frac{T_1}{T_2}s_{12}
 \ee
 where $T_1$ and $T_2$ are defined by
 \beqn
 T_1:&=&4m(2m-1)(m-1)C^m_{\frac{1}{2}}(-\frac{1}{b^2})^{m-1}+
 2(2m-1)(2m^2-2m-1)b^2a_{2m}\\
 &&+4(m+1)b^4a_{2m+2},\\
 T_2:&=&m(2m-1)^2b^3\big[a_{2m}-C^m_{\frac{1}{2}}(-\frac{1}{b^2})^m\big].
 \eeqn
 Now plug (\ref{yy76}), (\ref{yy83}), (\ref{yy84}) into
(\ref{yy49}), and then we obtain

 \be
  \frac{\phi-s\phi'
+(b^2-s^2)\phi''}{ s \phi + (b^2-s^2)\phi'} = \frac{ 2m-1}{s (1-k
s^2)},\label{y10****}
 \ee
  where $k$ is a constant. By the same
argument, we obtain (\ref{y10}).  This gives  Theorem
\ref{th3}(iv).

\section{Projectively flat $(\alpha,\beta)$-metrics}\label{sec7}
In this section, we characterize two-dimensional
$(\alpha,\beta)$-metrics (might be singular) which are
projectively flat. We have the following theorem.

\begin{thm}\label{th02}
  Let $F=\alpha \phi(s)$, $s=\beta/\alpha$, be an
  $(\alpha,\beta)$-metric on an open subset $U\subset R^2$ with $\phi(0)=1$. Suppose that
    $\beta$ is not parallel with respect to $\alpha$ and $F$ is not
of Randers type.  Then $F$ is  projectively flat in $U$ with
$G^i=P(x,y)y^i$ if and only if   $F$ lies in one of the following
four classes:
 \ben
  \item[{\rm (i)}] $\phi(s)$ and $\beta$ satisfy (\ref{yg3}) and
  (\ref{yg4}), and the spray coefficients $G^i_{\alpha}$ of $\alpha$
  satisfy
   \be\label{yg9}
   G^i_{\alpha}=\rho y^i-\tau (k_1\alpha^2+k_2\beta^2)b^i.
   \ee
 In this case, the projective factor $P$ is given by
  \be\label{yg10}
  P=\rho+\tau\alpha\Big\{\big[1+(k_1+k_3)s^2+k_2s^4\big]\frac{\phi'}{\phi}-(k_1+k_2s^2)s\Big\}.
  \ee

 \item[{\rm (ii)}] $\phi(s)$ and $\beta$ satisfy (\ref{y5}) and
  (\ref{y6}), and
  \be\label{g2}
   G^i_{\alpha}=\rho y^i-\tau \big\{(2c-k)\alpha^2+(c+k)k\beta^2\big\}b^i
   +\frac{(2c-k)(1-kb^2)\alpha^2+ck\beta^2}{1-(c+k)b^2}s^i.
   \ee
   In this case, the projective factor $P$ is given by
   \be\label{g02}
   P=\rho-\frac{4c^2s^3}{1+(c-k)s^2}\tau\alpha+\frac{\sigma_1}{\sigma_2}
   s_0,
   \ee
   where
   \beqn
   \sigma_1:&=&\big[(k-2c)(c-k)kb^2+4c^2-3kc+k^2\big]s^2
   +(2c-k)(1-kb^2),\\
   \sigma_2:&=&\big[1-(c+k)b^2\big]\big[1+(c-k)s^2\big].
  \eeqn

 \item[{\rm (iii)}] $\phi(s)$ and $\beta$ satisfy (\ref{y8}) and
  (\ref{y9}), and
 \be\label{g3}
  G^i_{\alpha}=\rho y^i-\frac{(2m-1)\alpha^2+k\beta^2}{2mb^2}s^i.
 \ee
 In this case, the projective factor $P$ is given by
 \be\label{g03}
  P=\rho-\frac{s(1-ks^2)\phi'+\big [  (2m-1) + k s^2\big ] \phi}{2mb^2\phi}s_0.
 \ee

\item[{\rm (iv)}] $\phi(s)$ and $\beta$ satisfy (\ref{y10}) and
  (\ref{y11}), and
 \be\label{g4}
  G^i_{\alpha}=\rho y^i-\frac{2(m-1)\alpha^2+ k \beta^2}{(2m-1)b^2}s^i.
 \ee
 In this case, the projective factor $P$ is given by
 \be\label{g04}
  P=\rho-\frac{s(1-kb^2)\phi'+\big[ 2(m-1) +ks^2\big ] \phi }{ (2m-1) b^2 \phi} s_0.
 \ee
  \een
In the above, $\rho= c_1(x)y^1+c_2(x)y^2$ is a 1-form.
\end{thm}

By Theorem \ref{th02}(ii), we can easily get the following
corollary.

\begin{cor}\label{cor52}
 Let $F=\alpha\pm \beta^2/\alpha$ be  two-dimensional
 $(\alpha,\beta)$-metric. Then $F$ is locally projectively flat if and only
 if  $\beta$ satisfies (\ref{y0006}) and $G^i_{\alpha}$ satisfy
 \be\label{y00060}
  G^i_{\alpha}=\rho y^i\mp 2\tau \alpha^2
   b^i-\frac{2\alpha^2}{b^2\mp 1}s^i.
 \ee
  In this case, the projective factor $P$ is given by
   \be\label{yg013}
   P=\rho-\frac{4s^3}{1\pm s^2}\tau\alpha-\frac{2(2s^2\pm
   1)}{(b^2\mp 1)(s^2\pm 1)}s_0.
   \ee
\end{cor}

To prove Theorem \ref{th02}, it follows from comparing Proposition
\ref{prop1} and Proposition \ref{prop2} that we only need to give
the expressions  (\ref{yg9})--(\ref{g04}) for each class in
Theorem \ref{th02}.

\subsection{The Spray Coefficients of $\alpha$}
In this subsection  we will show the expressions of the spray
coefficients $G^i_{\alpha}$ for each class in Theorem \ref{th02}.
Note that by $\widetilde{G}^i_{jk}=\frac{\pa^2G_{\alpha}^i}{\pa
y^j\pa y^k}$, the spray $G^i_{\alpha}$ of $\alpha$ can be
expressed as
 $$G^i_{\alpha}=\frac{1}{2}\widetilde{G}^i_{jk}y^jy^k.$$

\bigskip
\noindent
 {\bf Case I:}
Suppose that $(r_{11},r_{22})\ne (0,0)$. It has been proved in
\cite{Shen1} that
 \be\label{y96}
  \widetilde{G}^1_{11}=2t_1-2\lambda b^3\tau,\ \
  \widetilde{G}^1_{22}=-2\mu
  b^3\tau, \ \ \widetilde{G}^1_{12}=t_2, \ \
  \widetilde{G}^2_{12}=t_1,
 \ee
where $t_1,t_2$ are numbers independent of $s$ and $\tau$ is given
by (\ref{yg4}) or (\ref{y6}).  By (\ref{y42}) and (\ref{y43}) we
can get $\widetilde{G}^2_{11}$ and $\widetilde{G}^2_{22}$.

 If $\beta$ is closed ($s_{12}=0$), then it follows from (\ref{y32}), (\ref{y96}) and the
 expressions of $\widetilde{G}^2_{11}$ and $\widetilde{G}^2_{22}$
  that (\ref{yg9}) holds, where we put $\rho$ as  $\rho=t_iy^i$. This case has been given by
 \cite{Shen1} in case of $n\ge 3$.

If $\beta$ is not closed, then by putting $k_1 = 2c-k$ and $k_3 =
-3c -k$ we get (\ref{g2}) from (\ref{y52}), (\ref{y96}) and the
 expressions of $\widetilde{G}^2_{11}$ and $\widetilde{G}^2_{22}$.

\bigskip
\noindent
 {\bf Case II:}
Suppose that $(r_{11},r_{22})= (0,0)$. Then by (\ref{g26}) we get
 \be\label{g99}
  \widetilde{G}^1_{22}=0,\ \
  \widetilde{G}^1_{11}=2\widetilde{G}^2_{12}=2t_1, \ \
  \widetilde{G}^2_{12}=t_1,\ \ \widetilde{G}^1_{12}=t_2.
 \ee

If $db\ne 0$, then we have shown in Section \ref{sec5} that
$a_4\ne -\frac{1}{2}a_2^2$. In this case,
 we obtain (\ref{g2}) with $\tau=0$.

If $db= 0$, then we get (\ref{y01}). If $a_{2m+1}\ne 0$ in
(\ref{y01}), then we get from (\ref{yy80}), (\ref{yy81}) that
 \be\label{g100}
 \widetilde{G}^2_{11}=-\frac{2m-1+kb^2}{mb}s_{12},\ \
 \widetilde{G}^2_{22}=2\widetilde{G}^1_{12}-\frac{2m-1}{mb}s_{12}.
 \ee
Then it follows from (\ref{g99}) and
 (\ref{g100}) that (\ref{g3}) holds.
If $a_{2m}\ne C^m_{\frac{1}{2}}(-\frac{1}{b^2})^m$ in (\ref{y01}),
 then we get from (\ref{yy83}), (\ref{yy84}) that
 \be\label{g101}
 \widetilde{G}^2_{11}=-\frac{2(2m^2-2m-k)}{m(2m-1)b}s_{12},\ \
 \widetilde{G}^2_{22}=2\widetilde{G}^1_{12}-\frac{4(m-1)}{(2m-1)b}s_{12}.
 \ee
Then it follows from (\ref{g99}) and
 (\ref{g101}) that (\ref{g4}) holds.

 \subsection{The Projective Factors}
In this subsection,  we are going to find the expression for the
projective factor for each class in Theorem \ref{th02}.

Actually, (\ref{yg10}) has been proved in \cite{Shen1}, since
$\beta$ is closed in Theorem \ref{th02}(i). So we only show the
expressions of $P$ in (\ref{g02}), (\ref{g03}) and (\ref{g04}). In
the left three classes, since $\beta$ may not be closed, it is not
easy to show the projective factors $P$ in the initial local
projective coordinate system (in such a coordinate system,
geodesics are straight lines). However, it is easy to be solved by
choosing another local projective coordinate system, and then
returning to the the initial local projective coordinate system,
just as that in \cite{Shen1}.

Fix an arbitrary point $x_o\in U\subset R^2$. By the above idea
and a suitable affine transformation, we may assume $(U,x^i)$ is a
local projective coordinate system satisfying that
$\alpha_{x_o}=\sqrt{(y^1)^2+(y^2)^2}$ and $\beta_{x_o}=by^1$. Then
at $x_o$ we have
 $$s^1=s_1=0,\ \ s^2=s_2=bs_{12}, \ \ s_0=bs_{12}y^2, \ \ b^1=b_1=b, \ \
 b^2=b_2=0.$$

Suppose (\ref{y5}), (\ref{y6}) and (\ref{g2}) hold in $U$. Then it
is easy to get $r_{00},s_{0}^i,Q,\Theta$, and $\Psi$.
  Plug them into (\ref{y20}), and then at $x_o$
  we see that $G^i=Py^i$, where $P$ is given by
  \be\label{y107}
  P=\rho+A_1\tau+A_2bs_{12}y^2,
  \ee
 where
 $$A_1:=-\frac{4 c^2(by^1)^3}{[ 1+ (c-k)b^2] (y^1)^2+(y^2)^2},$$
 $$A_2:=\frac{[ 1+(2c-k)b^2] \big[ (2c-k)-k(c-k) b^2 \big](y^1)^2+(2c-k)(1-kb^2)(y^2)^2}
 {[1-(c+k)b^2] \big\{ [ 1+(c-k)b^2] (y^1)^2+(y^2)^2\big\}}.$$
 By using
   $$bs_{12}y^2=s_0,\ \ (y^1)^2+(y^2)^2=\alpha^2,\ \
  by^1=\beta,\ \ \frac{\beta}{\alpha}=s,$$
 we can transform (\ref{y107}) as
  (\ref{g02}). It is a direct computation, so the details are
  omitted. Since $x_o$ is arbitrarily chosen, (\ref{g02}) holds in
  $U$.

The left proofs are similar. So the details are omitted.

We have found the projective factor for each class in Theorem
\ref{th02}. This also gives a proof to the inverse of Theorem
\ref{th02}.

\section{Singular classes in Theorem \ref{th3} and  \ref{th02}}\label{sec6}

In this section, we will firstly discuss the local structures of
the singular classes in Theorem \ref{th3}(iii) and (iv), and then
construct some examples for Theorem \ref{th02}(iii) and (iv).

\

 Since every two-dimensional Riemann metric is locally conformally
 flat, we may put
  \be\label{y109}
  \alpha=e^{\sigma(x)}\sqrt{(y^1)^2+(y^2)^2},
  \ee
 where $x=(x^1,x^2)$. Since (\ref{y13}) is equivalent to
 $b^2=||\beta||^2_{\alpha}=constant$ in two-dimensional case (see a simple proof in \cite{LS2}),
 (\ref{y13}) holds if and only if $\beta$ is in the form
  \be\label{y110}
 \beta=\frac{be^{\sigma(x)}\big[\xi(x)y^1+\eta(x)y^2\big]}{\sqrt{\xi(x)^2+\eta(x)^2}},
  \ee
where $b=||\beta||_{\alpha}=constant$. Thus (\ref{y109}) and
(\ref{y110}) give all the local solutions of (\ref{y13}).
 If we put $\alpha$ and $\beta$ in the forms (\ref{y109}) and
(\ref{y110}), then we have (\ref{y13}), that is,
 $$
 r_{ij}=2\tilde{\tau}(b^2a_{ij}-b_ib_j)-\frac{1}{b^2}(b_is_j+b_js_i),$$
where $\tilde{\tau}$ is given by
 \be\label{y111}
 2\tilde{\tau}=\frac{(\xi^2+\eta^2)(\xi\sigma_1+\eta\sigma_2)-\xi\eta\eta_1+\xi^2\eta_2
 +\eta^2\xi_1-\xi\eta\xi_2}{be^{\sigma}(\xi^2+\eta^2)^{\frac{3}{2}}},
 \ee
where $\sigma_1:=\pa \sigma/\pa x^1$, $\sigma_2:=\pa \sigma/\pa
x^2$, etc.  Further, $\beta$ is not closed if and only if
 \be\label{y112}
(\xi^2+\eta^2)(\xi\sigma_2-\eta\sigma_1)-\xi^2\eta_1-\xi\eta\eta_2
 +\xi\eta\xi_1+\eta^2\xi_2\ne 0.
 \ee
In particular, if we put $\xi=x^2,\eta=-x^1$ and
$\sigma=c\big[(x^1)^2+(x^2)^2\big]$, where $c\ne 0$ is a constant,
then we have  $\tilde{\tau}=0$ by (\ref{y111}) and (\ref{y112}) holds
($\beta$ is not closed).

\begin{prop}
Define a two-dimensional $(\alpha,\beta)$-metric $F$  on $R^2$ by
 $$F=\frac{b\alpha^2+k\beta^2}{\sqrt{b^2\alpha^2-\beta^2}},$$
 where $k,b$ are constants with $k\not = -1/b$. Then $F$ is a Douglas metric if and only if $\alpha$ and $\beta$
 can be locally defined by (\ref{y109}) and (\ref{y110}), where
 $\xi,\eta$ and $\sigma$ are some scalar functions on $R^2$. There are
 many choices for $\xi,\eta$ and $\sigma$ such that $\beta$ is not
 closed.
\end{prop}

\begin{prop}
Let $F=\alpha \phi(s)$, $s=\beta/\alpha$  be a two-dimensional
  $(\alpha,\beta)$-metric, where $\phi(0)=1$. Let $\phi(s)$ be given by  (\ref{y8}) or (\ref{y10})(not given by (\ref{y13})).
   Then $F$ is a Douglas metric if and only if $\alpha$ and $\beta$
 can be locally defined by (\ref{y109}) and (\ref{y110}), where
 $\xi,\eta$ and $\sigma$ are some scalar functions satisfying $\tilde{\tau}=0$ in (\ref{y111}). There are
 many choices for $\xi,\eta$ and $\sigma$ such that $\beta$ is not
 closed.
\end{prop}

\

Next we construct some singular examples for Theorem
\ref{th02}(iii) and (iv)  which are  projectively flat. One can
directly verify the following two examples.

\begin{ex}\label{ex01}
Let $F=\alpha \phi(s)$, $s=\beta/\alpha$, be two-dimensional
  $(\alpha,\beta)$-metric, where $\phi(0)=1$. Let $\phi(s)$ be given
   by (\ref{y8}) with $k=0$, and define $\alpha$ and $\beta$ by
  (\ref{y109}) and (\ref{y110}), where
 $$\xi=x^2,\ \ \eta=-x^1,\ \
 \sigma=(m-\frac{1}{2})\ln\big[(x^1)^2+(x^2)^2\big].$$
  Then $F$ is projectively flat
  with $\beta$ being not closed.
\end{ex}
\begin{ex}\label{ex02}
Let $F=\alpha \phi(s)$, $s=\beta/\alpha$, be two-dimensional
  $(\alpha,\beta)$-metric, where $\phi(0)=1$. Let $\phi(s)$ be given by
   (\ref{y10}) with $k=0$, and define $\alpha$ and $\beta$ by
  (\ref{y109}) and (\ref{y110}), where
  $$\xi=x^2,\ \ \eta=-x^1,\ \ \sigma=(m-1)\ln\big[(x^1)^2+(x^2)^2\big].$$
   Then $F$ is projectively flat
  with $\beta$ being not closed.
\end{ex}

It might be also an interesting problem to show the local
structures of the two classes of Theorem \ref{th02}(iii) and (iv).
This problem is still open.

\section{Proof of Theorem \ref{th2}}

In the following proof, our idea is to choose a suitable
deformation on $\alpha$ and $\beta$ such that we can obtain a
conformal form on a Riemannian space. Then using the result in
\cite{Y}, we can complete our proof.

Define a new Riemann metric $\widetilde{\alpha}$ and a 1-form
$\widetilde{\beta}$ by
 \be\label{yg83}
\widetilde{\alpha}:=\sqrt{\xi\alpha^2+\eta\beta^2}, \ \
\widetilde{\beta}:=\beta,
 \ee
where
 $$
 \xi:=\frac{(1\mp b^2)^3}{(1\pm 2b^2)^{3/2}}, \ \
 \eta:=\frac{9}{8b^2}\Big\{(1\pm 2b^2)^{3/2}-\frac{1\mp
 2b^2+4b^4}{(1\pm 2b^2)^{3/2}}\Big\}.
 $$
 Since
$F=\alpha\pm\beta^2/\alpha$ is a Douglas metric, we have
(\ref{y0006}). By (\ref{yg83}) and (\ref{y0006}), a direct
computation gives
 \be\label{yg85}
 \widetilde{r}_{ij}=\frac{2\tau (1\mp b^2)^2}{(1\pm
 2b^2)^{5/2}}\widetilde{a}_{ij}.
 \ee
So $\widetilde{\beta}=\beta$ is a conformal 1-form with respect to
$\widetilde{\alpha}$.

Since $\widetilde{\alpha}$ is a two-dimensional Riemann metric, we
can express $\widetilde{\alpha}$ locally as
 \be\label{yg86}
\widetilde{\alpha}:=e^{\sigma}\sqrt{(y^1)^2+(y^2)^2},
 \ee
where $\sigma=\sigma(x)$ is a scalar function. We can obtain the
local expression of  $\widetilde{\beta}=\beta$ by (\ref{yg85}) and
(\ref{yg86}) (see \cite{Y}). By the result in \cite{Y}, we have
 \be\label{yg87}
 \widetilde{\beta}=\widetilde{b}_1y^1+\widetilde{b}_2y^2=e^{2\sigma}(uy^1+vy^2),
 \ee
where $u=u(x),v=v(x)$ are a pair of scalar functions such that
 $$f(z)=u+iv, \ \ z=x^1+ix^2$$
 is a complex analytic function.

 We can express $\sigma$ using $u,v,b^2$ by computing the quantity
 $||\beta||^2_{\widetilde{\alpha}}$. Firstly, by (\ref{yg86}) and
(\ref{yg87}) we get
 \be\label{yg88}
||\beta||^2_{\widetilde{\alpha}}=e^{2\sigma}(u^2+v^2).
 \ee
On the other hand, by the definition  of $\widetilde{\alpha}$ in
(\ref{yg83}), the inverse $\widetilde{a}^{ij}$ of
$\widetilde{a}_{ij}$ is given by
 $$
\widetilde{a}^{ij}=\frac{1}{\xi}\Big(a^{ij}-\frac{\eta
b^ib^j}{\xi+\eta b^2}\Big).
 $$
Now plug $\xi$ and $\eta$ into the above, and we obtain
 \be\label{yg89}
||\beta||^2_{\widetilde{\alpha}}=\widetilde{a}^{ij}b_ib_j=\frac{b^2}{(1\pm
2b^2)^{\frac{3}{2}}}.
 \ee
Thus by (\ref{yg88}) and (\ref{yg89}) we get
 \be\label{ycw93}
 e^{2\sigma}=\frac{1}{u^2+v^2}\frac{b^2}{(1\pm
 2b^2)^{\frac{3}{2}}}.
 \ee

Finally, by plugging (\ref{ycw93}), (\ref{yg86}) and (\ref{yg87})
into (\ref{yg83}), we easily get $\alpha$ and $\beta$ given by
(\ref{yc1}) and (\ref{yc2}) respectively, where we define
$B:=b^2$.
   \qed

\

The Riemann metric $\alpha$ expressed in (\ref{yc1}) is generally
not in the conformally flat form. We show another  example which
is expressed in a different form. Put
 \beqn
 \alpha&=&e^{\sigma}\sqrt{(y^1)^2+(y^2)^2}, \ \ \
\beta=e^{\sigma}(\xi y^1+\eta y^2),\\
\sigma&=&-\frac{3}{2}\ln\big[1\pm 2(x^1)^2\pm 2(x^2)^2\big]+c,\\
\xi&=&x^2,\ \ \eta=-x^1,
 \eeqn
 where $c$ is a constant. Then $\beta$ is not closed and satisfies
(\ref{y0006}) with $\tau=0$. Therefore, $F=\alpha\pm
\beta^2/\alpha$ is a Douglas metric.

\section{Proof of Theorem  \ref{th001}}

This proof is similar as that in Theorem \ref{th2}. In the proof
of Theorem \ref{th2}, we make a deformation on $\alpha$ and keep
$\beta$ unchanged. For the proof of Theorem  \ref{th001} in the
following, we will give a deformation on $\beta$ but keep $\alpha$
unchanged.

Define a  Riemannian metric $\widetilde{\alpha}$ and 1-form
$\widetilde{\beta}$ by
 \be\label{yg90}
\widetilde{\alpha}:=\alpha, \ \
\widetilde{\beta}:=\frac{\beta}{c},
 \ee
where $c=c(b^2)$ is defined in (\ref{yg79}), where we define
$B:=b^2$. By (\ref{yg4}) and a direct computation we can obtain
 \be\label{yg91}
 \widetilde{b}_{i|j}=\frac{2\tau (1+k_1b^2)}{c}\widetilde{a}_{ij}=\frac{2\tau (1+k_1b^2)}{c}a_{ij}.
 \ee
So $\widetilde{\beta}$ is a closed 1-form conformal  with respect
to $\alpha$.

Now we  express $\alpha$ locally as
 \be\label{yg92}
\alpha:=e^{\sigma}\sqrt{(y^1)^2+(y^2)^2},
 \ee
where $\sigma=\sigma(x)$ is a scalar function.  Then by the result
in \cite{Y}, we have
 \be\label{yg93}
 \widetilde{\beta}=\widetilde{b}_1y^1+\widetilde{b}_2y^2=e^{2\sigma}(uy^1+vy^2),
 \ee
where $u=u(x),v=v(x)$ are a pair of scalar functions such that
 $$f(z)=u+iv, \ \ z=x^1+ix^2$$
 is a complex analytic function.

 By the property of $u,v$ and the fact that $\widetilde{\beta}$ given in (\ref{yg93}) is
 closed, it is easily seen that $u,v,\sigma$ satisfy the PDEs (\ref{yg81}).

 Now we determine $\sigma$ in terms of the triple $(B,u,v)$, where $B:=b^2$.
 Firstly by (\ref{yg90}) and then by (\ref{yg92}) and
(\ref{yg93}) we get
 \be\label{yg94}
||\widetilde{\beta}||^2_{\alpha}=b^2c^{-2},\ \
||\widetilde{\beta}||^2_{\alpha}=e^{2\sigma}(u^2+v^2).
 \ee
Therefore, by (\ref{yg94}) we get (\ref{sigma}).

Finally, we can easily get $\alpha$ and $\beta$ given by
(\ref{yg79}) from (\ref{sigma}), (\ref{yg90}), (\ref{yg92}) and
(\ref{yg93}).         \qed

\vspace{0.6cm}

\noindent Guojun Yang \\
Department of Mathematics \\
Sichuan University \\
Chengdu 610064, P. R. China \\
 ygjsl2000@yahoo.com.cn

\end{document}